\documentclass[letterpaper,12pt,leqno,oneside]{article}

\usepackage{float}
\usepackage{caption}
\usepackage{enumitem}
\usepackage{color}
\usepackage{ifsym}
\usepackage{bm}
\usepackage{exscale}
\usepackage{amsmath}
\usepackage{amsfonts}
\usepackage{wasysym}
\usepackage{stmaryrd}
\usepackage{amscd}
\usepackage{graphicx}
\usepackage{pagecolor,lipsum}
\usepackage{xcolor}
\usepackage{marvosym}
\usepackage{textcomp}
\usepackage{amsxtra}
\usepackage{amssymb}
\usepackage{theorem}
\usepackage{enumitem}
\usepackage{mathtools}
\usepackage{esint}
\usepackage[hidelinks]{hyperref}
\usepackage[final]{epsfig}
\usepackage{eqnarray}
\usepackage{hyperref}
\newtheorem{proposition}{Proposition}[subsection]
\newtheorem{definition}[proposition]{Definition}

{\theorembodyfont{\rmfamily}}
\newtheorem{theorem}[proposition]{Theorem}
\newtheorem{corollary}[proposition]{Corollary}

{\theorembodyfont{\rmfamily}}

\newfont{\abc}{cmtt10 scaled 1200}

\def\ve{\varepsilon}

\def\ra{\rightarrow}
\def\cs{\symbol{35}}
\def\p{\partial}

\def\mm{\mbox}
\def\v{= \emptyset}

\def\bp{\langle A \rangle}

\def\si{$\mathcal{S}$}
\def\sima{\mathcal{S}}

\def\R{\mathbb{R}}

\def\T{\mathbb{T}}
\def\Z{\mathbb{Z}}

\def\ve{\varepsilon}

\def\ve{\varepsilon}

\def\ra{\rightarrow}
\def\cs{\symbol{35}}
\def\p{\partial}
\def\bp{\langle A \rangle}

\begin{document}

\vspace*{-3cm}
\begin{center}\Large{\bf{The Secret Hyperbolic Life of Positive Scalar Curvature}}\\
\bigskip
\large{\bf{Joachim Lohkamp}}\\
\end{center}
{\small\noindent Mathematisches Institut, Universit\"at M\"unster, Einsteinstrasse 62, Germany\\
{\emph{e-mail: j.lohkamp@uni-muenster.de}}}\\

{\small \textbf{Abstract} \, This survey introduces to the \emph{hyperbolic unfolding} correspondence that links the geometric analysis of minimal hypersurfaces with that of Gromov hyperbolic spaces. Problems caused from hypersurface singularities oftentimes become solvable on associated Gromov hyperbolic spaces. Applied to scalar curvature geometry this yields \emph{smoothing schemes} that eliminate such singularities. We explain the two main lines of such smoothings: a top-down singular analysis using iterated blow-ups of the given hypersurface and bottom-up smoothings following the tree of blow-ups backwards.}

{\footnotesize  {\center \tableofcontents}}

\setcounter{section}{1}
\renewcommand{\thesubsection}{\thesection}
\subsection{Introduction} \label{introduction}

The story began in the 70s when Hawking \cite[Ch.~4]{H}  and  Schoen and Yau \cite{SY1} discovered the first cases of a remarkable $scal >0$-heredity one may state, in dimensions $n \ge 2$, as follows: \\

\emph{Given a compact manifold $M^{n+1}$ with $scal >0$, any area minimizing hypersurface $(H^n,g_H)  \subset (M^{n+1},g_M)$ also has  $scal >0$, after conformal deformations of $g_H$.}\\

We indicate how to derive this powerful yet simple result. The first and the second variation of $Area(H)$ of a hypersurface $H \subset M$ by a $f \cdot \nu$, where $\nu$ is the outward normal vector field of $H$, $f \in C^\infty(H,\R)$ with $supp \: f \subset reg(H)$ (= the set of regular points of $H$), are given by
\begin{align*}
    Area'(f)  &=  \int_H tr A_H (z) \cdot f(z) \: dVol, \\
    Area''(f) &=  \int_H |\nabla_H f|^2  + \left( (tr A_H)^2 -  |A|^2 - Ric(\nu,\nu) \right) \cdot f^2 \: dVol
\end{align*}
where $tr A_H$ is the mean curvature of $H$, $|A|^2$ is the sum of the squares of principal curvatures of $H$ and $Ric$ the Ricci tensor of $M$.
Since $H$ is supposed to be area minimizing, we have $Area'(f) = 0 $ and $Area''(f) \ge 0 $. That is, $tr A_H = 0$ and this gives
 {\small \begin{equation}
         \label{var2}  \quad Area''(f) = \int_{H}|\nabla_H f|^2  -  \left( |A|^2 + Ric(\nu,\nu)  \right) \cdot f^2 \: dVol \ge 0\end{equation}}
 Now we recall the Gau\ss--Codazzi equations for hypersurfaces
\[\label{gc} \small |A|^2 + Ric(\nu,\nu)  =1/2 \cdot \left(| A |^2 + scal_M - scal_H +(tr A_H)^2 \right)\] where $scal_H$ and $scal_M$ denote the scalar
curvature of $H$ and $M$. Since $tr A_H = 0$ we may rewrite (\ref{var2}) as follows (here we assume $n \ge 3$, the computations slightly deviate when $n=2$): \begin{multline} \label{3} \int_H | \nabla f |^2 + \frac{n-2}{4 (n-1)} \cdot scal_H \cdot f^2 d A \\
\ge \int_H \frac{n}{2 (n-1)} \cdot  |  \nabla f |^2 + \frac{n-
2}{4 (n-1)}\cdot \left( | A |^2 + scal_M \right) \cdot f^2 d A.
\end{multline}
The left hand side of $(\ref{3})$ is the variational integral for the first eigenvalue $\lambda_1$ of the conformal Laplacian $L_H$ on $H$. That is, when $scal_M > 0$, inequality $(\ref{3})$ shows that $\lambda_1 > 0$. Following Kazdan and Warner [KW], the transformation law for scalar curvature under conformal transformation by the first eigenfunction $f_1>0$ can be used to show that
\[ \label{4} scal(f_1^{4/(n-2)} \cdot g_H) \cdot f_1^{\frac{n+2}{n-2}} =
L_H(f_1) = \lambda_1 \cdot f_1 > 0.\]
In other words, we observe a reproduction of the scalar curvature constraint on $M$ on the lower dimensional space $H$, also called a scalar curvature \emph{splitting factor} of $M$. This suggests an appealing and versatile strategy to study scalar curvature from an inductive dimensional splitting until one reaches a space already understood. \\
The problem is that only in low dimensions $n \le 7$ these hypersurfaces are smooth submanifolds. The ultimate goal of this paper is to explain how to use hyperbolic geometry to get such smooth splitting factors in arbitrary dimensions. The basic reference are the papers \cite{L1}--\cite{L5}. The program is illustrated in the following diagram.
\begin{figure}[htbp]
\centering
\includegraphics[width=1 \textwidth]{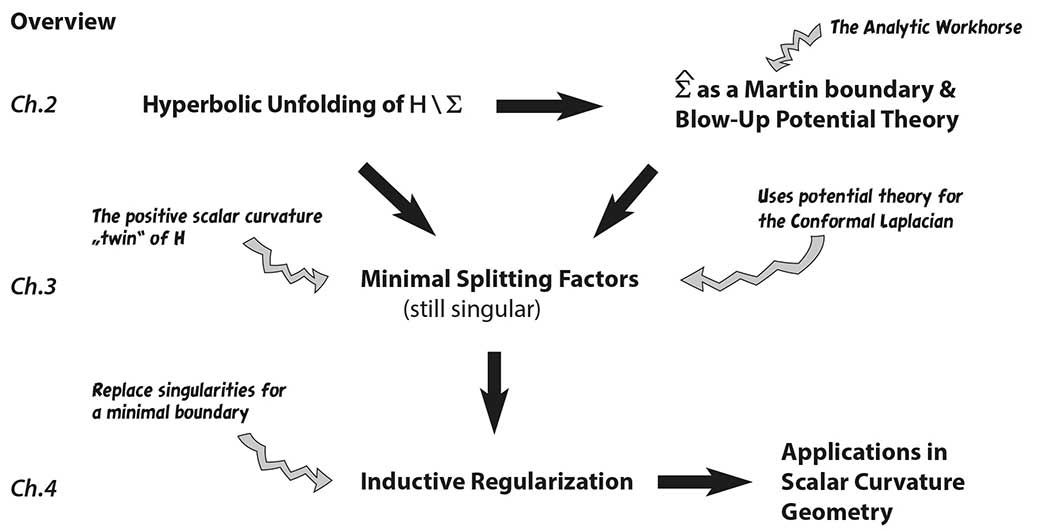}
\caption{The hyperbolic unfolding gives use means to control the elliptic analysis on $H$ but it also reveals some new geometric structures like canonical Semmes families on $H$.}
\end{figure}
Hyperbolic geometry is extensively used right from the beginning, in chapter 2. The smoothings occupy chapters 3 and 4.\\
As a general remark, the theory we present naturally extends to considerably larger classes of almost minimizers. Some of them do not even result from variational principles. Nevertheless, for the most part we are talking about area minimizers to simplify the exposition.

\setcounter{section}{2}
\renewcommand{\thesubsection}{\thesection}
\subsection{Singularities as Regular Boundaries} \label{hyunf}

The problem we encounter in dimensions $\ge 8$ is how to control $L_H$ near the singular set of $H$. The critical dimension $8$ comes from classical regularity theory. It says that an area minimizing hypersurface $H^n$ that is smooth outside a potentially complicated \textbf{singular set} $\Sigma_H \subset M^{n+1}$ of codimension $ \ge 8$. To overcome these subtleties we reinterpret $\Sigma_H$ as a \textbf{boundary} of $H \setminus \Sigma_H$ and show that, \textbf{relative} to $H \setminus \Sigma_H$, the singular set $\Sigma_H$ becomes \textbf{regular}, in a sense we explain below.\\
We make use of some pieces of geometric measure theory but we neither employ any structural details of $\Sigma$ nor the full codimension estimate but only that $\Sigma$ has codimension $\ge 3$. \\

\textbf{Motivation} \, To motivate and properly describe the adequate notion of regularity of $\Sigma $ we consider the case of Euclidean domains.  The validity of the \textbf{boundary Harnack inequality} on a Euclidean domain $G \subset \R^n$, in (\ref{gf2}) below, is a litmus test for a reasonable elliptic analysis on $G$. It asserts that there are constants $A(G)$, $C(G)>1$ such that for any  $p\in \p G$, small $R>0$ and any two harmonic functions $u$, $v>0$ on $B_{A \cdot R}(p) \cap G$ which vanish along $B_{A \cdot R}(p) \cap \p G$,
\begin{equation}\label{gf2}
u(x)/v(x) \le C \cdot  u(y)/v(y) \mm{ for all } x,\, y \in B_R(p) \cap G.
\end{equation}
Remarkably, there is a purely geometric way to characterize such domains. The validity of  (\ref{gf2}) is essentially \textbf{equivalent} \cite{Ai} to the regularity condition that $G$ is a \textbf{uniform domain}.

\begin{definition}\label{ud}  A Euclidean domain $G \subset \R^n$  is called   \textbf{uniform} if there is some $c\ge 1$ such that any two points $p,q \in G$ can be joined by a \textbf{uniform curve}. This is a rectifiable curve $\gamma_{p,q}: [a,b] \ra G$ going from $p$ to $q$ such that:
\begin{itemize}[leftmargin=*]
  \item \emph{\textbf{Quasi-geodesic:}} \, $l(\gamma_{p,q})\le c\cdot d(p,q)$.
  \item \emph{\textbf{Twisted double cones:}} \, For any $z \in \gamma_{p,q}$ let $l_{min}(\gamma_{p,q}(z))$ be the minimum of the lengths of the two subcurves of $\gamma_{p,q}$ from $p$ to $z$ and from $z$ to $q$. Then
  \[
  l_{min}(\gamma_{p,q}(z)) \le c \cdot dist(z,\p X).
  \]
\end{itemize}
\end{definition}
One may describe uniformity as a quantitative and scaling invariant form of connectivity.\\

\textbf{\si-Structures} \, Turning to $H \setminus \Sigma_H$ relative to its boundary $\Sigma_H$, we observe that $G \subset \R^n$ is flat until we reach $\p G$ whereas $H \setminus \Sigma_H$ degenerates towards $\Sigma_H$ with diverging second fundamental form. The isoperimetric inequality for the area minimizer $H$ allows us to compensate the additional twist from the divergence of second fundamental form and to establish the even stronger  \textbf{\si-uniformity} of $H$. To make this precise we introduce a measure $\bp_H$ for this degeneration.
An assignment $H \mapsto \bp_H$ of a non-negative, measurable function to any connected area minimizing hypersurface $H$ is called an \textbf{\si-transform} provided
\begin{itemize}[leftmargin=*]
\item $\bp_H$ is \textbf{naturally} assigned to $H$, in other words, the assignment commutes with the convergence of sequences of underlying area minimizers.
    \item $\bp_H \ge |A_H|$ with  $\bp_H \equiv 0$, if $H \subset M$ is totally geodesic, otherwise, $\bp_H$ is  strictly positive.
    \item When $H$ is not totally geodesic, the  \textbf{\si-distance} $\delta_{\bp}:=1/\bp$  is well-defined and it is
        $L_{\bp}$-Lipschitz regular, for some constant $L_{\bp}=L(\bp,n)>0$:
\begin{equation}\label{delip}
|\delta_{\bp}(p)- \delta_{\bp}(q)|   \le L_{\bp} \cdot d(p,q), \mm{ for } p,q \in  H \setminus \Sigma.
\end{equation}
\end{itemize}
The \si-distance is a measure for the distance to singular and highly curved portions of $H$ that takes also the curvature into account. Constructions of \si-transforms are given in \cite{L1} where we merge  $g_H$ and $A_H$ in a suitable way. The geometric main application is the following result.

\begin{theorem}\label{thm2} There exists some $c>0$ such that $H \setminus \Sigma$ is an  \textbf{\si-uniform space}. This means that any pair $p,q \in H \setminus \Sigma$ can be joined by an \textbf{\si-uniform curve} in $H \setminus \Sigma$, i.e., a rectifiable curve $\gamma_{p,q}: [a,b] \ra H \setminus \Sigma$ with $\gamma_{p,q}(a)=p$, $\gamma_{p,q}(b)=q$ and so that:
  \begin{itemize}[leftmargin=*]
    \item \emph{\textbf{Quasi-geodesic:}} \, $l(\gamma)  \le c \cdot  d(p,q).$
    \item \emph{\textbf{Twisted double \si-cones:}} \, $l_{min}(\gamma_{p,q}(z)) \le c \cdot \delta_{\bp}(z)$ for any $z \in \gamma_{p,q}$.
  \end{itemize}
\end{theorem}
\begin{figure}[h]
\centering
\includegraphics[width=1\textwidth]{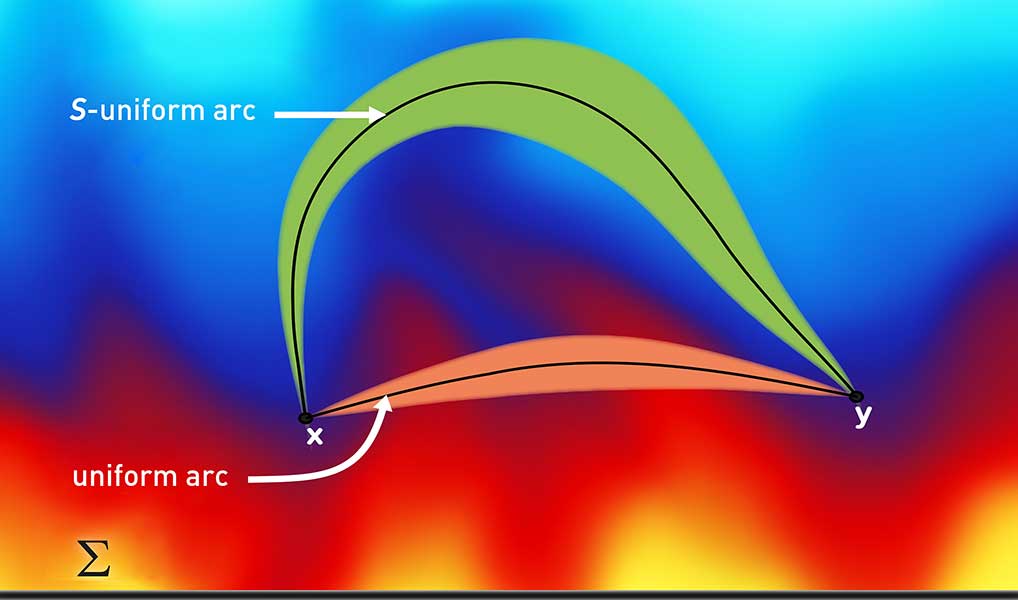}
\caption{Schematic view of $H$, where $\Sigma$ is just the bottom line. The cold colors indicate low
and the hot color high curvature of $H$. \si-uniform curves are also sensitive to the underlying curvature.}
\end{figure}
\textbf{Some Remarks} \, In this theory the totally geodesic hypersurfaces  play the  r\^ole of the trivial case. They are always smooth submanifolds \cite[Corollary A.6]{L1} and in the non-compact case they are just Euclidean hyperplanes. In this case, many results in this paper are either obvious or they degenerate to conventions.\\
The Lipschitz regular \si-distance $\delta_{\bp}$ admits a Whitney type $C^\infty$-smoothing $\delta_{\bp^*}$ which satisfies (S1)--(S3) and is quasi-natural in the sense that $c_1 \cdot \delta_{\bp}(x) \le \delta_{\bp^*}(x)  \le c_2 \cdot \delta_{\bp}(x)$, for some constants $c_1$, $c_2>0$, cf.\ \cite[Proposition B.3]{L1}.\\
We note in passing that in older papers, now superseded by \cite{L1}--\cite{L3}, we used the term \emph{skin transform} for the subclass of \si-transforms that includes a so-called Hardy inequality, now called \emph{Hardy \si-transforms}, cf.~\cite{L3}. The renaming had no deeper reasons. The Hardy inequality is only needed in applications but not in the basic theory. Secondly, originally the level sets of $\bp_H$, the \emph{$|A|$-skins}, served important technical purposes now covered from functional relations we have for $\bp_H$.

\subsubsection{Hyperbolic Unfoldings} \label{ung}

Remarkably, the uniformity of $G$ is also \textbf{equivalent} to the purely geometric condition that the  \textbf{quasi-hyperbolic metric}, which we get from conformally deforming the Euclidean metric $g_{Eucl}$ to $dist(\cdot,\p G)^{-2} \cdot g_{Eucl}$, is Gromov hyperbolic,  has bounded geometry and its Euclidean boundary is homeomorphic to the Gromov boundary, cf.~\cite{BHK}. To make this plausible, we first observe that a (twisted) cone conformally deformed by $dist(\cdot,\p G)^{-2}$ roughly looks like a piece of a hyperbolic space (a generalized Poincar\'{e} metric) and the scaling invariance of the uniformity conditions globalizes this to the whole space.\\
We will see that an area minimizer $H$ with singular set $\Sigma$ (= boundary of $H \setminus \Sigma$) has a similar hyperbolic nature. This is more challenging since $H \setminus \Sigma$ degenerates while we approach  $\Sigma$.  In zeroth order this can be  compared with the extension of the Riemann uniformization from complex domains to arbitrary Riemann surfaces.  \\

For starters we recall the notions of Gromov hyperbolicity and of bounded geometry.  It has no local impact but strong consequences for the geometry near infinity. We mention \cite{BH} as a general reference.

\begin{definition}\label{del}
A locally compact geodesic metric space $X$ is \textbf{Gromov hyperbolic} if its geodesic triangles are $\mathbf{\delta}$\textbf{-thin} for some $\delta=\delta_X >0$. That is, each point on the edge of any geodesic triangle is within $\delta$-distance of one of the other two edges.\\
Two rays in $X$ are \emph{equivalent} if they have finite Hausdorff distance. The set $\p_G X$ of equivalence classes $[\gamma]$ of geodesic rays from a fixed base point $p \in X$ is called the \textbf{Gromov boundary} of $X$. This definition of $\p_G X$ is independent of $p$. The space $\overline{X}_G = X \cup \p_G X$ admits a natural topology that makes it a compact metrizable space. It is called the \textbf{Gromov compactification} of $X$.
\end{definition}

\begin{definition}\label{bog}
A Riemannian manifold $M$ has $(\varrho,\ell)$-\textbf{bounded geometry} if there exist constants $\varrho={\varrho_M}>0$ and $\ell=\ell_M \ge 1$ for $M$ such that for each ball $B_{\varrho}(p) \subset M$ there is a smooth $\ell$-bi-Lipschitz chart $\phi_p$ onto an open set $U_p\subset\R^n$ with its Euclidean metric.\\
\end{definition}

\textbf{The \si-metric on $H \setminus \Sigma$} \,  In general the quasi-hyperbolic metric $dist(\cdot,\Sigma_H)^{-2} \cdot g_H$ is not well-behaved on $H \setminus \Sigma$. It neither has bounded geometry nor does it have a good blow-up behavior, that is, the corresponding  quasi-hyperbolic metric on tangent cones of points in $\Sigma_H \subset H$ does not approximate that on $H \setminus \Sigma$. The key point about \si-uniformity of $H$ is that it implies all the desirable properties for the so-called \textbf{\si-metric} $d_{\bp}=d_{\bp_H}$  defined by
\begin{equation}\label{sim}
d_{\bp}(x,y) := \inf \Bigl  \{\int_\gamma  \bp \, \, \Big| \, \gamma   \subset  H \setminus \Sigma\mbox{ rectifiable curve joining }  x \mbox{ and } y  \Bigr \}
\end{equation}
for $x$, $y\in H \setminus \Sigma$. The metric $d_{\bp}$ is even well-defined for \emph{smooth} $H$ where $\Sigma\v$. Alternatively, the \si-metric can be written $\bp^2 \cdot g_H$, but this is not a regular Riemannian metric since $\bp$ is merely a locally Lipschitz function. However, one may use the Whitney type $C^\infty$-smoothing $\delta_{\bp^*}$ if one needs a smooth version. Now we can formulate the following hyperbolization result, cf.\ \cite[Theorem 1.11, Proposition 3.10 and Theorem 1.13]{L1}.

\begin{theorem}
The \si-metric $d_{\bp}$ has the following properties:
\begin{itemize}[leftmargin=*]
  \item The metric space $(H \setminus \Sigma, d_{\bp})$ and its \textbf{quasi-isometric Whitney smoothing}, i.e., the smooth Riemannian manifold $(H \setminus \Sigma, d_{\bp^*}) = (H \setminus \Sigma, 1/\delta_{\bp^*}^2 \cdot g_H)$, are \textbf{complete Gromov hyperbolic spaces} with \textbf{bounded geometry}.
  \item $d_{\bp}$ is \textbf{natural}, that is, the assignment $H\mapsto d_{\bp_H}$ commutes with compact convergence of regular domains of the underlying area minimizers. The typical example is the blow-up around singular points and the resulting  tangent cone approximations.
  \item For any \emph{singular} $H$ the identity map on $H \setminus \Sigma$ extends to \textbf{homeomorphisms}
  \begin{align*}
  \widehat{H}&\cong\overline{(H \setminus \Sigma,d_{\bp})}_G \cong \overline{(H \setminus \Sigma,d_{\bp^*})}_G\,\mm{ and }\\
  \widehat{\Sigma} &\cong\p_G(H \setminus \Sigma,d_{\bp}) \cong \p_G(H \setminus \Sigma,d_{\bp^*}),
  \end{align*}
  where for $X=(H \setminus \Sigma,d_{\bp})$ or $(H \setminus \Sigma,d_{\bp^*})$, $\overline{X}_G$ and $\p_G(X)$ denote the Gromov compactification and the Gromov boundary, respectively.
\end{itemize}
The spaces $(H \setminus \Sigma, d_{\bp})$ and $(H \setminus \Sigma, d_{\bp^*})$ are conformally equivalent to the original space $(H \setminus \Sigma, g_H)$. We refer to both these spaces as \textbf{hyperbolic unfoldings} of $(H \setminus \Sigma, g_H)$.
\end{theorem}

Here, $\widehat{H}$ and $\widehat{\Sigma}$ denote the one-point compactifications of $H$ and $\Sigma$ in the non-compact case of Euclidean hypersurfaces with the extra condition that for $H \subset \R^{n+1}$ we always add the point $\infty$
to $\Sigma$, even for compact $\Sigma$.

\subsubsection{Reduction to Hyperbolic Geodesics} \label{rtg}

Hyperbolic unfoldings reveal that the \textbf{potential theory} of elliptic operators $L$, like the conformal Laplacian or the Jacobi field operator, on area minimizers \textbf{conforms to the underlying geometry}. Intuitively, the evolution of the minimal Green's function $\widetilde{G}$ of $\widetilde{L}=\delta_{\bp^*}^2\cdot L$ on $(H \setminus \Sigma, d_{\bp^*})$ concentrates along hyperbolic geodesics and this it largely controls the global analysis of $\widetilde{L}$. (The minimal Green's function $\widetilde{G}$ of $\widetilde{L}$,
$G(x,y)>0$ is a smallest function on $M \times M$, singular on $\{(x,x)\:|\:x \in M\}$, satisfying the equation $\widetilde{L}\,G(\cdot,y)=\delta_y$, where $\delta_y$ is the Dirac function in $y$.)\\
This focussing/squeezing effect of hyperbolic geodesics is reflected in so-called $\boldsymbol{3G}$\textbf{-inequalities}, due to the triple appearance of $G$ in one inequality, one may interpret as follows:\\

\emph{Let $x, y, z  \in (H \setminus \Sigma, d_{\bp^*})$ such that $y$ lies on a hyperbolic geodesic connecting $x$ and $z$. Then, up to universal constants, there are as many
``Brownian particles'' travelling directly from $x$ to $z$, measured by $\widetilde{G}(x,z)$, as there are particles travelling from $x$ to $y$, measured by $\widetilde{G}(x,y)$, and then from $y$ to $z$, measured by $\widetilde{G}(y,z)$. That is, we have $\widetilde{G}(x,z) \approx\widetilde{G}(x,y) \cdot \widetilde{G}(y,z)$.}\\

The potential theory on hyperbolic manifolds of bounded geometry naturally extends to hyperbolic graphs of bounded valence, where the constant $\delta$ measures the deviation of the graph from a tree. In turn, we can approximate $(H \setminus \Sigma, d_{\bp^*})$, along with its potential theory, by such graphs. This formalizes the intuition of a reduction to hyperbolic geodesics.\\

The admissible operators on $(H \setminus \Sigma, d_{\bp^*})$ to make this work are the \textbf{adapted weakly coercive} operators. These are the linear second order elliptic operators $L$ which are \textbf{uniformly elliptic} on $(H \setminus \Sigma, d_{\bp^*})$ and \textbf{weakly coercive}. That is, there is a $u >0$ with $L \, u \ge \ve \cdot u$ for  some $\ve  >0.$  For an exposition of this potential theory on Gromov hyperbolic spaces, which is largely due to Ancona \cite{A1,A2}, see \cite{KL} for an exposition.\\
An interesting aspect of this theory is that the geometric and the analytic conditions work hand in hand: the pair of conditions \emph{bounded geometry} $\leftrightarrow$ \emph{uniformly elliptic} is employed for most of the basic estimates and relations and then these results are critically improved from a tandem use of \emph{Gromov hyperbolicity} $\leftrightarrow$ \emph{weak coercivity}.\\

The point about hyperbolic unfoldings is that the transparent potential theory we have on $(H \setminus \Sigma, d_{\bp^*})$ transfers one-by-one to corresponding classes of operators on the original area minimizer. This way we get the same type of potential theoretic results we have for uniform Euclidean domains on the \si-uniform spaces $H \setminus \Sigma$ with the difference that we have literally outsourced all the analysis to the hyperbolic unfolding as our  workbench. This is the \textbf{unfolding corrrespondence} and the relevant operators $L$ on $(H \setminus \Sigma, g_H)$ are called \emph{\si-adapted} provided their counterpart
$\delta_{\bp^*}^2\cdot L$ on $ (H \setminus \Sigma, d_{\bp^*})$ is adapted weakly coercive. More explicitly, we set
\begin{definition}\label{sadp}
An elliptic operator $L$ is \textbf{\si-adapted} when $-L(u) = \sum_{i,j}  a_{ij} \cdot \frac{\p^2 u}{\p x_i \p x_j} + \sum_i b_i \cdot \frac{\p u}{\p x_i} + c \cdot u$ so that for some $k \ge 1$ and suitable charts
\begin{equation}\label{1}
k^{-1} \cdot\sum_i \xi_i^2 \le \sum_{i,j} a_{ij}(p) \cdot \xi_i \xi_j \le k \cdot \sum_i \xi_i^2,
\end{equation}
\begin{multline}\label{2}
\delta^{\beta}_{\bp}(p) \cdot  |a_{ij}|_{C^\beta(B_{\Theta(p)}(p))} \le k, \delta_{\bp}(p) \cdot |b_i|_{L^\infty(B_{\Theta(p)}(p))} \le k \mm{ and }\\
\delta^2_{\bp}(p) \cdot |c|_{L^\infty(B_{\Theta(p)}(p))} \le k, \mm{ for } \Theta(p)=c/\bp(p), \mm{ for some } c(H)>0.
\end{multline}
and there is a $u >0$ with $L \, u \ge \ve \cdot \bp^2 \cdot u\mm{ for  some } \ve  >0.$
\end{definition}

A frequently considered type of elliptic problems is that of eigenvalues, typically for \textbf{symmetric operators}. For them it is possible and useful to bring the weak coercivity condition into a variational form.  Then the weak coercivity condition  is equivalent to the existence of some positive constant $\tau>0$  such that the Hardy type inequality
\begin{equation}\label{hadi0}
\int_H  f  \cdot  L f  \,  dV \, \ge \, \tau \cdot \int_H \bp^2\cdot f^2 dV \mm{ holds for any } f \in C^\infty_0(H \setminus \Sigma).
\end{equation}
That is, $L$ has an $\bp$-weighted \textbf{positive first eigenvalue} $\lambda^{\bp}_1(L)>0$, in this singular case commonly called the\textbf{ principal eigenvalue}. We also mention that in the singular case, in sharp contrast to the smooth case, we have positive $\bp$-weighted eigenfunctions for any $\lambda < \lambda^{\bp}_1(L)$.\\

The boundary Harnack inequality we get for \si-adapted $L$ along the boundary $\widehat{\Sigma}$ of $(H\setminus \Sigma, g_H)$ differs in two ways from that in  the case of uniform domains equipped with the Laplacian.
\begin{itemize}[leftmargin=*]
  \item In place of the balls of (\ref{gf2}), we choose hyperbolic halfspaces $\mathcal{N}^\delta_i(z)$ in $(H \setminus \Sigma, d_{\bp^*})$, with completions $\mathbf{N}^\delta_i(z)$, uniformly contracting to $z \in \widehat{ \Sigma}$ for $i \ra \infty$. One may compare this with the Poincar\'{e} disc model on $B_1(0) \subset \R^{2}$ where the $B_\rho(p) \cap B_1(0)$, $\rho \in (0,1)$,  for flat discs $B_\rho(p)$, $p \in \p D$, become hyperbolic halfspaces.
\item Solutions $u >0$ of  $L \, f=0$ will usually diverge to infinity when we approach $\Sigma$. The  generalization of the vanishing boundary data is that of solutions of \textbf{minimal growth} towards $\Sigma_H$ when compared to  other solutions $v>0$.
\end{itemize}

\begin{theorem} There exists a constant $C(H,L) >1$ such that for any $z \in \widehat{\Sigma}$ and any  two solutions
 $u$, $v >0$ of $L\, f= 0$ on $H \setminus \Sigma$ with minimal growth along $\mathbf{N}^\delta_i(z)\cap \widehat{\Sigma}$, we have
\begin{equation}\label{fhepq1}
u(x)/v(x) \le C \cdot  u(y)/v(y) \mm{ \emph{for all} }x,\, y \in \mathcal{N}^\delta_{i+1}(z).
\end{equation}
\end{theorem}

\textbf{Examples} \, To relate the general setup of \si-adapted operators to some geometrically relevant examples, we upgrade the \si-transform $\bp_H$ to a \textbf{Hardy \si-transform}, cf.~\cite{L3}.  This is an \si-transform that additionally satisfies the following Hardy type inequality:   for any $f \in C^\infty(H \setminus \Sigma,\R)$
        compactly supported in $H \setminus \Sigma$ we have
\[\int_H|\nabla f|^2  + |A_H|^2 \cdot f^2 dA \ge \tau \cdot \int_H \bp_H^2\cdot f^2 dA, \mm{ for some } \tau = \tau(\bp,H) \in (0,1).\]

Let $H^n \subset M^{n+1}$ be a singular area minimizer. Then we have for any Hardy \si-transform:
\begin{itemize}[leftmargin=*]
\item For $scal_M \ge 0$ we have $\lambda^{\bp}_1(L_H) >0$ for the \textbf{conformal Laplacian} $L_H$.
\item We have $\lambda^{\bp}_1(J_H) \ge 0$ for the \textbf{Jacobi field operator} $J_H:=-\Delta_H - |A|^2-Ric_M(\nu,\nu)$.
\end{itemize}
On area minimizers, we can turn any operator that satisfies (\ref{1}) and  (\ref{2}) into an \si-adapted operator $L$ provided $\lambda^{\bp}_1(L) > -\infty$ since $L - \lambda \cdot \bp^2 \cdot Id$ becomes \si-adapted for $\lambda < \lambda^{\bp}_1(L)$.\\

\subsubsection{Applications in Geometric Analysis} \label{aot}

A basic application of the boundary Harnack inequality is a transparent Martin theory.
\begin{definition}\label{mb}
Let $X$ be a non-compact Riemannian manifold and $L$ be a linear second order elliptic operator on $X$ with a minimal Green's function $G: X\times X \ra (0,\infty]$. We choose a base point $p$ and consider the space $S$ of sequences $s=\{p_n\}$ in $X$, $n \ge 1$, such that
\begin{itemize}[leftmargin=*]
  \item $s$ has no accumulation points in $X$.
  \item $K(x,p_n):= G(x,p_n)/G(p,p_n) \ra K_s(x)$ compactly to some function $K_s$ on $X$ as $n \ra \infty$.
\end{itemize}
The \textbf{Martin boundary} $\p_M (X,L)$ is the quotient of $S$ modulo the following relation on $S$: $s \sim s^*$ if and only if $K_s \equiv K_{s^*}$. Moreover, we define the \textbf{Martin kernel} $k(x;y)$ on $X \times \p_M (X,L)$ by $k(x;y):= K_s(x)$, for some sequence $s$ representing $y \in \p_M (X,L)$. As for the Gromov boundary, these definitions do not depend on the choice of the base point $p$.
\end{definition}

The (metrizable) \textbf{Martin topology} on $\overline{X}_M:= X \cup \p_M (X,L)$ is defined through the convergence of the function $K(x,p_n)$. $\overline{X}_M$ and $\p_M (X,L)$  turn out to be compact. $\overline{X}_M$ is called the \textbf{Martin compactification} of $(X,L)$.   A point in $S_L(X)$ is called \textbf{extremal} if it cannot be written as a non-trivial convex combination of other points of $S_L(X)$. The subset $\p^0_M (X,L) \subset \p_M (X,L)$ of extremal functions belonging to $\p_M (X,L)$ is called the \textbf{minimal Martin boundary}. \\
The point about $\p^0_M (X,L)$ is that for any solution  $v>0$  of $L\,f=0$ we have
a \textbf{unique} finite Radon measure $\mu=\mu(v)$ on $\p^0_M (X,L)$ so that
\begin{equation}\label{mint2}
v(x)  =\int_{\p^0_M (X,L)} k(x;y) \, d \mu(y) .
\end{equation}
The problem with this \textbf{Martin integral} is that a proper understanding of $\p^0_M (X,L)$ can be very difficult. Different from classical contour formulas, $\p_M (X,L)$  and $\p^0_M (X,L)$ may strongly depend on $L$ and they usually differ from intrinsically defined topological boundaries of $X$.  Remarkably,  these problems disappear for \si-adapted operators on area minimizers.

\begin{theorem}\label{mbhs}
For any  \si-adapted operator $L$ on $H \setminus \Sigma$ we have
\begin{itemize}[leftmargin=*]
  \item the identity map on $H \setminus \Sigma$ extends to a homeomorphism between $\widehat{H}$ and the Martin compactification $\overline{(H \setminus\Sigma)}_M$.
  \item all Martin boundary points are minimal:  $\p^0_M (H \setminus \Sigma,L) \equiv \p_M(H \setminus \Sigma,L)$.
\end{itemize} Thus,  $\widehat{\Sigma}$ and the minimal Martin boundary $\p^0_M (H \setminus \Sigma,L)$ are homeomorphic.
\end{theorem}

\textbf{Quantitative Results} \, The unfolding correspondence is not restricted to the transition between $L$ on $(H \setminus \Sigma, g_H)$ and $ \delta_{\bp^*}^2\cdot L$ on $(H \setminus \Sigma, d_{\bp^*})$. Inspired from the Doob transform ín stochastic analysis we define
\begin{equation}\label{lsdef}
L^{\sima} \phi: = \delta^{(n+2)/2}_{\bp^*} \cdot  L ({\delta^{-(n-2)/2}_{\bp^*}} \cdot \phi) \mm{ for sufficiently regular } \phi \mm{ on } (H \setminus \Sigma, d_{\bp^*})
\end{equation}
We call $L^{\sima}$ the \textbf{\si-Doob Transform}  of $L$.
The minimal Green's functions $G$ of an \si-adapted $L$ on $(H \setminus \Sigma, d_H)$ and $G^{\sima}$ of the adapted weakly coercive $L^{\sima}$ on $(H \setminus \Sigma, d_{\bp^*})$ satisfy
\begin{equation}\label{g1}
G^{\sima}(x,y) =  \delta_{\bp^*}^{(n-2)/2}(x)  \cdot   \delta_{\bp^*}^{(n-2)/2}(y) \cdot G(x,y), \mm{ for } x \neq y \in H \setminus \Sigma.
\end{equation}
There are constants $\beta(H),\alpha(H)>0$  so that (typically applied along hyperbolic geodesics)
\begin{multline}\label{g2}
 G^{\sima}(x,y)\le \beta  \cdot \exp(-\alpha \cdot d_{\bp^*}(x,y)),\mm{ for } x,y  \in H \setminus \Sigma \mm{ and } \\d_{\bp^*}(x,y)> 2 \cdot \varrho_{(H \setminus \Sigma, d_{\bp^*})}.
\end{multline}
In the case of the conformal Laplacian, the combination of \eqref{g1} and \eqref{g2} with the axioms for $\bp$ yields upper radius and distance estimates after conformal deformations
by $G$ and, via boundary Harnack inequality, for arbitrary solutions of minimal growth.\\

\textbf{Minimal Growth under Blow-Ups} \,  In the case of a singular minimal cone $C$  the Martin theory shows that there is exactly one solution $u>0$ with minimal growth towards $\Sigma_C$, up to multiples. This means that $\mu_u$ is the Dirac measure in $\infty \in \widehat{\Sigma}_C$. Now we assume that $L$ reproduces under scalings, $L_C$ and $J_C$ are examples. This implies a separation of variables: $u = \psi_C(\omega) \cdot r^{\alpha_C}, \, (\omega,r) \in  \p B_1(0)  \cap C \setminus \Sigma_C \times \R^{>0}$, for some function $\psi_C$ on $\p B_1(0)  \cap C \setminus \Sigma_C$ and $\alpha_C <0$.\\
When $u$ solves $L \, v = 0$ with minimal growth along $B \cap \widehat{\Sigma}$ for a ball $B \subset H$ around some $p \in \Sigma$ and we
consider any solution $v$ on any tangent cone $C$ of $H$ in $p$ induced by $u$ under blow-up, we find that $v$ has minimal growth towards all of $\Sigma_C$.
$u$ asymptotically looks like  $\psi_C(\omega) \cdot r^{\alpha_C}$. The outcome is a \textbf{top-down} to \textbf{bottom-up} asymptotic analysis towards $p$.\\

\textbf{$\bullet$ Top-Down} \, Given a solution $u>0$ on $H \setminus \Sigma_H$ of  minimal growth near some $p \in \Sigma_H$, choose some tangent cone $C$ and consider the (uniquely determined) induced solution of minimal growth towards $\Sigma_C \subset C$ on $C \setminus \Sigma_C$. Blow-up in a point of $\Sigma_C \setminus \{0\}$ and iterate this process, at most $\dim H -7$ times, until we reach a product cone $\R^m \times C^{n-m} $, for some $C^{n-m} \subset \R^{n-m+1}$ singular only in $0$. We think of this as a \textbf{terminal node} in a blow-up tree with \textbf{root} $H$. The uniqueness shows that the induced minimal solutions we end up with are $\R^{n-m+1}$-translation symmetric and amenable to an explicit description.\\

\textbf{$\bullet$ Bottom-Up} \, Starting from terminal nodes $\R^m \times C^{n-m} $, for some $C^{n-m} \subset \R^{n-m+1}$ singular only in $0$, we transfer our understanding of the induced solutions on the tangent cones stepwise backwards to $H$ in the blow-up tree.  We get an asymptotic portrait of the solution $u$ on  $H \setminus \Sigma$ near $p$ from the family of Martin theories on $H \setminus \Sigma$ and on its (iterated) tangent cones.\\

\setcounter{section}{3}
\renewcommand{\thesubsection}{\thesection}
\subsection{Minimal Splitting Factors - Top-Down Analysis}

For a compact singular area minimizer $H^n$  in some $scal>0$-manifold $M^{n+1}$ we recall that $\lambda^{\bp}_H>0$.
We apply our theory to get a nicely controllable conformal but still singular $scal>0$-geometry on $H^n$, the \textbf{minimal splitting factor} geometry.
One may think of it as some kind of  \textbf{$\boldsymbol{scal>0}$-twin} of $H^n$ and we will explain why.

\subsubsection{Minimal Splitting Factors} \label{mm}

To this end, we can neither use $L_H$ nor $L_{H,\lambda^{\bp}_H } = L_H - \lambda^{\bp}_H \cdot \bp^2 \cdot Id$. When we deform $H$, more precisely $H \setminus \Sigma$, by solutions of $L_H \phi =0$ we get a scalar flat metric. One may feel tempted to consider solutions of $L_{H,\lambda^{\bp}_H} \phi =0$, but this operator is no more \si-adapted simply because it obviously has a vanishing principal eigenvalue and many of the results do not apply.
Instead we choose (super)solutions of  $L_{H,\lambda} \phi= L_H \phi - \lambda \cdot \bp^2 \cdot \phi =0$ for
\begin{equation}\label{lam}
0 < \lambda < \lambda^{\bp}_H.
\end{equation}

Due to the locally Lipschitz regular coefficients of $L_{H,\lambda} = L_H - \lambda \cdot \bp^2 \cdot Id$, solutions  of $L_{H,\lambda} \, \phi=0$ are $C^{2,\alpha}$-regular, for any $\alpha \in (0,1)$. This suggests the following regularity assumptions. For  $\lambda < \lambda^{\bp}_H$, let  $\Phi>0$ be a $C^{2,\alpha}$-supersolution  of $L_{H,\lambda} \phi =0$ on $H\setminus\Sigma_{H}$ so that:
\begin{itemize}[leftmargin=*]
  \item  For compact area minimizers $H$: $\Phi$  is a solution in a neighborhood of $\Sigma$   and it has \textbf{minimal growth} towards $\Sigma$.
  \item For non-compact Euclidean area minimizers $H$:  $\Phi$ is a solution on $H\setminus \Sigma_{H}$ with \textbf{minimal growth} towards $\Sigma$.
\end{itemize}

\begin{theorem}  The metric completion $(\widehat{H \setminus \Sigma}, \widehat{d_{\sima}}(\Phi))$ of $(H \setminus \Sigma,\Phi^{4/(n-2)} \cdot g_H)$ is \textbf{homeomorphic} to $(H,d_H)$. Thus, we can write it as $(H,d_{\sima}(\Phi))$.
The \textbf{Hausdorff dimension} of $\Sigma$ relative to $(H,d_{\sima}(\Phi))$ is $\le n-7$.
\end{theorem}
Here $(H,d_H)$ is the
metric space that results from the embedding $H^n \subset M^{n+1}$. We write simply $d_{\sima}$ when the specific choice of $\Phi$  is not needed or already known from the context.\\

\textbf{Some Details} \, To get the homeomorphism, we use the \si-Doob transform (\ref{lsdef}) and the relation (\ref{g1}) to transfer the basic upper growth estimate (\ref{g2}) to the minimal Green's function of $L_{H,\lambda}$ on $(H \setminus \Sigma, g_H)$. It also applies to $\Phi$ by the boundary Harnack inequality. This shows that $(\widehat{H \setminus \Sigma}, \widehat{d_{\sima}})$ is not larger than $(H,d_H)$.  This works as soon as $\lambda < \lambda^{\bp}_H$. Remarkably, we need $\lambda >0$ to also get lower estimates showing that it is not smaller and, hence, to establish the homeomorphism from the Bombieri--Giusti Harnack inequality \cite{BG}.\\
One might expect that the Hausdorff dimension of $\Sigma$ relative to $(H,d_{\sima})$ must be the same as for area minimizers. However, the Hausdorff dimension is not a topological but only a bi-Lipschitz invariant while the identity map from $(H, d_H)$ to $(H, d_{\sima})$ is not Lipschitz continuous. H\"{o}lder continuous homeomorphisms can increase the dimension of a subset of dimension $a \in (0,n)$ to any value $b \in (a,n)$, cf.~\cite{Bi}. We actually need the minimal growth of $\Phi$  and essential parts of the hyperbolic unfolding theory to control the Hausdorff dimension.\\

\textbf{Metric Measure Spaces} \,  We augment $(H, d_{\sima}(\Phi_H))$ to a metric measure space. To this end we show that there is a canonical extension of $\Phi^{2 \cdot n/(n-2)}\cdot \mu_H$ on $H  \setminus \Sigma$ to a measure $\mu_{\sima}$ on $H$, where $\mu_H$ is the $n$-dimensional Hausdorff measure on $(H^n,g_H) \subset (M^{n+1},g_M)$.  In fact, this extension $\mu_{\sima}$ is a \textbf{Borel measure} on $(H,d_{\sima})$, cf.~\cite[pp.\,62--64]{HKST}.

\begin{definition} We call  $(H,d_{\sima})$ a \textbf{minimal spitting factor}  of its ambient space $M$  and we define the \textbf{minimal factor measure} $\mu_{\sima}$ by $\mu_{\sima}(E):=\int_{E \setminus \Sigma_H} \Phi^{2 \cdot n/(n-2)}\cdot d\mu_H,$ for any Borel set $E \subset H$.
\end{definition}
The small Hausdorff dimension of $\Sigma \subset (H, d_{\sima}(\Phi_H))$ ensures that $\mu_{\sima}$ is an \textbf{outer regular measure} and it is still sufficient to define $\mu^{n-1}_{\sima}$ for hypersurfaces  within $(H, d_{\sima}(\Phi_H))$. We get $d\mu^{n-1}_{\sima}$ from extending $\Phi^{2 \cdot (n-1)/(n-2)}\cdot d\mu^{n-1}_H$ on $H \setminus \Sigma$, where $d\mu_H^{n-1}$ is the hypersurface element  on $(H \setminus \Sigma, g_H)$.

\begin{theorem} We consider $(H, d_{\sima}(\Phi_H))$, some $p \in \Sigma_H$ and some tangent cone $C$ in $p$. Then we get the following \textbf{blow-up invariance:}\\
Any sequence  $(H, \tau_i \cdot d_{\sima}(\Phi_H))$, scaled by a sequence $\tau_i \ra \infty$, $i \ra \infty$, around $p$, subconverges and the limit of any converging subsequence is $(C, d_{\sima}(\Phi_C))$ for some tangent cone $C$.   $(C, d_{\sima}(\Phi_C))$ is invariant under scaling around $0 \in C$, that is, it is again a cone.
 \end{theorem}

This is the geometric counterpart of the  \textbf{top-down} blow-up analysis of minimal growth solutions we discussed in the previous chapter. The result means that $(H,d_{\sima},\mu_{\sima})$ has a simple asymptotic geometry near $\Sigma$ and it admits inductive tangent cone reductions similar to area minimizers. Note that the principal eigenvalues of $H$ and $C$ usually differ, with $\lambda^{\bp}_{H} < \lambda^{\bp}_{C}$. That is, we would encounter non-principal eigenvalues in the blow-up analysis even if we started from $\lambda^{\bp}_{H}$.\\

\subsubsection{Stable Isoperimetry - Hyperbolic Geodesics Again} \label{csf}

A vital feature of minimal spitting factors is that they  still satisfy an isoperimetric inequality.

\begin{theorem} There are some constants $\gamma(H)>0,$ $\gamma(H^*)>0$ so that for every ball $B_\rho$ and any open set $U \subset H$ with compact closure and rectifiable boundary $\p U$
\begin{equation}\label{iifin}
 \mu_{\sima}(U)^{(n-1)/n} \le \gamma \cdot \mu^{n-1}_{\sima}(\p U),
\end{equation}
\begin{equation}\label{ii2in}
\min  \{ \mu_{\sima}(B_{\rho} \cap U),  \mu_{\sima} (B_{\rho} \setminus U)\}^{(n-1)/n} \le \gamma^* \cdot \mu^{n-1}_{\sima}(B_{\rho} \cap \p U).
\end{equation}
\end{theorem}

The proof has two steps. We first show that $(H,d_{\sima},\mu_{\sima})$ is doubling and has a volume decay property of order $n$.
Both follow from the stronger \textbf{Ahlfors regularity} estimate that for some $a(H), b(H)>0$,
\begin{equation}\label{volds}
 a \cdot r^n \le  Vol(B_r(q),d_{\sima}) \le b \cdot r^n.
\end{equation}
The trick is to check that $Vol(B_1(q),d_{\sima})$ for Euclidean hypersurfaces can be bounded in terms of constants that only depend on the dimension. This uses the same techniques as the proof of the homeomorphism theorem above. Scalings to other radii commute with scaling of the original area minimizing metric and readily give these growth rates.\\

The main step is the proof of the Poincar\'{e} inequality on $(H,d_{\sima},\mu_{\sima})$. It says that there is a constant $C_0=C_0(H,\Phi) >0$, so that when $B \subset H$ is an open ball, $u:B \ra \R$ is an $L^1$-function on $H$ that is $C^1$ on $H \setminus \Sigma$. Then we have, setting $|\nabla u| \equiv 0$ on $\Sigma$,
\begin{multline}\label{poinm}
\fint_B |u-u_B| \,  d \mu_{\sima} \le C_0 \cdot \fint_B |\nabla u| \, d \mu_{\sima},\\
 \mm{ where } u_B:=\fint_B u  \,  d \mu_{\sima} := \frac1{\mu_{\sima}(B)}\int_B u \, d \mu_{\sima}.
\end{multline}
The volume decay property of order $n$ then improves this Poincar\'{e} inequality to the Sobolev inequality with exponent $n$ and, from this, we get the isoperimetric inequalities for $(H,d_{\sima},\mu_{\sima})$. A broad reference is \cite{HKST}.\\

\textbf{Semmes Families} \, To derive the Poincar\'{e} inequality  on rather general metric measure spaces, Semmes has decompiled the classical proof of Poincar\'{e} inequality on $\R^n$ \cite{Se,He,HKST}. In an important step, we encounter uniformly distributed families of curves linking any two given points.
The abstracted concept is that of \emph{thick families of curves}, also called \emph{Semmes families}, satisfying the two conditions (i) and (ii)  below. For reasonable metric spaces, the presence of Semmes families implies the validity of a Poincar\'{e} inequality. On area minimizers, we can
use hyperbolic unfoldings to define  \textbf{canonical Semmes families} on  $(H, d_H)$ with the interesting feature that they are still Semmes families relative to $(H, d_{\sima})$.

\begin{proposition}  There is  a $C=C(H)>0$
and for any two $p,q \in H$ a family $\Gamma_{p,q}$ of rectifiable curves $\gamma: I_\gamma \ra H$, $I_\gamma \subset \R$,  joining $p$ and $q$, so that:
\begin{enumerate}
 \item For any $\gamma \in \Gamma_{p,q}$: $l(\gamma|_{[s,t]}) < C \cdot d(\gamma(s),\gamma(t))$, for  $s,t \in I_\gamma$.
 \item Each family $\Gamma_{p,q}$ carries a probability measure $\sigma_{p,q}$ so that for any Borel set $A \subset X$, the assignment $\gamma \mapsto l(\gamma \cap A)$  is $\sigma$-measurable with
 {\small \begin{equation}\label{tcu1}
 \int_{\Gamma_{p,q}} l(\gamma \cap A) \, d \sigma(\gamma) \le  C \cdot \int_{A_{C,p,q}} \left(\frac{d(p,z)}{\mu(B_{d(p,z)}(p))} + \frac{d(q,z)}{\mu(B_{d(q,z)}(q))}\right) d \mu(z)
 \end{equation}}
 for $A_{C,p,q}:=(B_{C \cdot d(p,q)}(p) \cup B_{C \cdot d(p,q)}(q))\cap A$.
\end{enumerate}
The family $\Gamma_{p,q}$ uniformly surrounds a central curve, its \textbf{core} $\gamma_{p,q}$. It is a hyperbolic geodesic linking $p$ and $q$ in the Gromov compactification of the hyperbolic unfolding.
 \end{proposition}

The first assertion, for the core, is part of the proof that $(H \setminus \Sigma, d_{\bp^*})$ is Gromov hyperbolic in \cite{L1}. It shows that each segment of  $\Gamma_{p,q}$ is an \si-uniform curve relative to $(H,d_H)$. For the other curves, this follows from the way we define $\Gamma_{p,q}$.\\

To define  $\Gamma_{p,q}$ and the probability measure $\sigma_{p,q}$, we start with  any two points $x, y \in \R^n$ and consider the hyperplane $L^{n-1}(x,y)$ orthogonal to the line segment $[x, y] \subset \R^n$ passing through the midpoint $m(x,y)$ of $[x, y]$. For $\rho=d(x,y)$, we consider a ball $B_r=B_r^{n-1}(m(x,y)) \subset L^{n-1}(x,y)$ of radius $r  \in (0,\rho]$ which we will choose later. For any $z\in B_r$, let $\gamma_{z}$ be the unit speed curve from $x$ to $y$ we get when we follow the line segments $[x, z]$ and
$[z,y]$. We define $\Gamma^{\R^n}_{x,y}=\Gamma^{\R^n}_{x,y}(r):=\{\gamma_{z}\,|\,z\in B_r\}$ and the point set $D_r(x,y)$ we get as a union of the curves.
We introduce the probability measure $\alpha^{\R^n}_{x,y}$ on $\Gamma^{\R^n}_{x,y}$ as follows: if $W\subset\Gamma^{\R^n}_{x,y},$
\begin{equation}\label{defsem}
\alpha^{\R^n}_{x,y}(W)  :=\mathcal{H}^{n-1}(\{z\in B_r^{n-1}\,|\, \gamma_{z}\in W\})/\mathcal{H}^{n-1}(B_r^{n-1}).
\end{equation}
For any Borel set $A\subset \R^{n}$, the function $\gamma\mapsto\ell(\gamma\cap A)$ on $\Gamma^{\R^n}_{x,y}$ is $\alpha^{\R^n}_{x,y}$-measurable.
Now assume that all points in $A$  are closer to $x$ than to $y$. For the distance between corresponding points on any two segments, we have
\begin{equation}\label{dist}
d(s \cdot z_1 + (1-s) \cdot x, s \cdot z_2 + (1-s) \cdot x) \le 2 \cdot r \cdot s, \mm{ for } s \in [0,1], z_1, z_2 \in B_r.
\end{equation}
The coarea formula gives the following inequality  for the annuli $A_{j} :=A \cap B(x,\ 2^{-j}d(x,\ y))\setminus B(x, 2^{-j-1}d(x, y))$,
$j \in \Z^{\ge 0}$:

\[\int_{B_r}\ell(\gamma_{z}\cap A)d\mathcal{H}^{n-1}(z)=\sum_{j=0}^{\infty}\ \int_{B_r}\ell(\gamma_{z}\cap A_{j})d\mathcal{H}^{n-1}(z) \le  \sum_{j=0}^{\infty} 2^{(j+2)(n-1)} \mu(A_{j})\]
\[\le 4^{n-1} \cdot \int_{A\cap B(x,d(x,y))}\frac{d(x,z)}{\mu(B_{d(x,z)}(x))} d \mu(z)\]

Thus we have a Semmes family on $\R^n$ and can now identify the line segment $[x, y] \subset \R^n$  with the hyperbolic geodesic $\gamma_{p,q}$ and use its \si-uniformity to transfer the Semmes family $\Gamma^{\R^n}_{x,y}$ to $(H, d_H)$.
In other words, what we really do is to define a \textbf{curve fibration} of the twisted double \si-cone surrounding $\gamma_{p,q}$ which we already had from the \si-uniformity of $H$. This is the desired canonical Semmes family on  $(H, d_H)$.\\
Once again, the estimates for minimal Green's functions along hyperbolic geodesics show that,
with other constants and probability measures, $\Gamma_{p,q}$  is still a  Semmes family for $(H,d_{\sima},\mu_{\sima})$ and thus we get the desired Poincar\'{e} inequality on $(H,d_{\sima},\mu_{\sima})$.\\

The  isoperimetric inequalities also control area minimizers within $(H,d_{\sima},\mu_{\sima})$ from the following type of standard consequences:
\begin{corollary}\label{horn} There
is a  $\rho_H>0$ so that for any $r \in (0,\rho_H)$ and any area minimizing boundary $L^{n-1}$ bounding some open set $L^+ \subset H$ in $(H,d_{\sima},\mu_{\sima})$:
\begin{equation}\label{estin}
 \kappa^-_n\cdot r^n  \le  \mu_{\sima}(L^+ \cap B_r(p)) \le  \kappa^+_n\cdot r^n,
\end{equation}
 where $\kappa^-_n, \kappa^+_n  >0$ denote constants depending only on the dimension $n$.
\end{corollary}

This volume growth estimate shows that there are no horn-shaped pieces of area minimizers in $(H,d_{\sima},\mu_{\sima})$ entering narrow hollow cylinders. In particular, such an area minimizer cannot \emph{stretch out} along $\Sigma$.
This is a crucial detail to improve the efficiency of smoothing techniques we discuss in the next chapter.

\setcounter{section}{4}
\renewcommand{\thesubsection}{\thesection}
\subsection{Smoothing Techniques - Bottom-up Constructions}

We describe how to deform the still singular minimal splitting factors $(H,d_{\sima})$ to regular $scal>0$-geometries on $H \setminus U$ with minimal boundary $\p U$, where is $U$ is a small neighborhood of $\Sigma_H$. This yields schemes of an inductive dimensional descent with a built-in partial regularization.\\
The point is that $\p U$ may also have singularities but of lower dimension and the minimality of $\p U$ matches the
intended inductive descent:\\  Area minimizers in $H \setminus U$  either coincide with $\p U$ or they are entirely supported in $H \setminus \overline{U}$.
Thus in each loop of the induction such area minimizers experience a smooth $scal>0$-environment and we eventually reach exclusively smooth geometries in dimensions $\le 7$.

\subsubsection{Surgeries Revisited} \label{c3}

We shed some new light on the well-known $scal >0$-preserving surgeries on a manifold $M$ one can do along any given submanifold $N \subset M$ of codimension $\ge 3$, \cite{GL} and \cite{SY2}.
We approach these surgeries in a way that represents the gluing boundaries, after removing tubes around $N$, as minimal surfaces. One may keep the intermediate bounded manifolds per se or glue some complementary piece to  them to get  closed $scal >0$-manifolds with altered topology.\\
This illustrates the smoothing results for minimal splitting factors we discuss below. In a quite similar manner, we remove tubes around the singularities $\Sigma \subset H$. Again,
we will get minimal boundaries and keep $scal >0$. One may use this manifold with minimal boundary directly or one can glue a smooth complementary part to it to get a new closed $scal >0$-manifold.

\begin{proposition} \label{cd3b}  Let $M^n$, $n \ge 3$, be a  closed  manifold and assume that the first eigenvalue $\lambda_1$ of the conformal Laplacian $L_M$ is positive. For any
submanifold $N^k \subset M^n$ of dimension $k \le n-3$ there are arbitrarily small neighborhoods $U$ of $N$,  such that $M \setminus U$ is conformal to a $scal>0$-manifold $X=X_U$ a (locally) area minimizing boundary $\p X$.
\end{proposition}

 We break the construction of the conformal deformation into two steps: In Step A we choose a canonical conformal deformation to reach a \emph{basic $scal>0$-metric} on $M$. In Step B we add a\emph{ secondary} conformal change of $M \setminus N$ using some positive function diverging to $+\infty$ while we approach $N$ to achieve the desired \emph{bending} effect towards $N$. When we evaluate these conformal deformations, we find two neighborhoods $U \subset V$ of $N$ with $\p V$ having \emph{positive} mean curvature and $\p U$ with \emph{negative} mean curvature\footnote{The sign convention is made so that $\p B_1(0) \subset \R^{n+1}$, viewed from $0 \in \R^{n+1}$,  has positive mean curvature.}.\\

We think of  $\p V$ as an \textbf{outer barrier} and $\p U$ as an \textbf{inner barrier}. They keep  area minimizers from escaping $V \setminus U$. In fact, the mean curvature constraints show that we may replace any $Area_{n-1}$-minimizing sequence of boundaries $T_i \subset V$ \emph{surrounding} $U$, that is $T_i = \p W_i$, for open  $W_i$ with $U \subset W_i \subset V$, so that for $i \ra \infty$:
\[Area_{n-1}(T_i) \ra \inf \{Area_{n-1}(T) \,|\, T = \p W\mm{ for  an open } W,  U \subset W \subset V \},\]
by another $Area_{n-1}$-minimizing  sequence $T^*_i$ with $Area_{n-1}(T^*_i) \le Area_{n-1}(T_i)$ and support outside a neighborhood of $\p V \cup \p U$. Then standard arguments from basic geometric measure theory \cite[1.20]{Gi} show that there is an area minimizer $T ^*= \p W^*$  for some open  $W^*$ with $U \subset W^* \subset V$. We set $X:=M \setminus W^*$ and have $\p X = T^*$.\\

\textbf{Step A}  \,  We can use the first eigenfunction u=$u_M>0$ of $L_M$ on $M$, the assumption $\lambda_1>0$ and the transformation law for scalar curvature under conformal transformation to see:
\begin{equation} \label{eig} scal(u^{4/(n-2)} \cdot g_M) \cdot u^{\frac{n+2}{n-2}} =
L_M(u) = \lambda_1 \cdot u > 0.\end{equation}
That is $ scal(u^{4/(n-2)} \cdot g_M) >0$. This observation is due to Kazdan and Warner \cite{KW}. $u^{4/(n-2)} \cdot g_M$ is our basic $scal>0$-metric on $M$.\\

\textbf{Step B}  \,  We  also choose a function $\psi>0$ on $M \setminus N$ with $L_M \psi =  0$ and
\begin{equation}\label{appp}
\psi=
    \begin{cases*}
 r^{2+k-n} + O( r^{3+k-n}), & when $k < n-3$ \\
r^{2+k-n} + O(log(r)^{-1}), & when $k = n-3$,
    \end{cases*}
  \end{equation}
where $r$ is the distance function to $N$. A construction of $\psi$ is given in \cite[Appendix]{SY2}. Roughly speaking,  one averages a Green's function for $L_M$ along $N$.   \\

Now we add $\delta \cdot \psi$ to $u$, for some small $\delta >0$, as a \textbf{secondary deformation}: we choose $g_\delta:=(\delta \cdot \psi + u)^{4/(n-2)} \cdot g_M$. The linearity of $L_M$ shows that $scal(g_\delta)>0$.\\
We denote the $r$-distance neighborhoods of $N$ relative to $u^{4/(n-2)} \cdot g_M$ by $V_r$. Then, for any small $\ve>0$, we have that $\p V_{\ve}$ is diffeomorphic to $N^k \times S^{n-k-1}$. Thus for small enough $\delta(\ve)>0$ we readily see that  $\p V_{\ve}$ is essentially determined from $u^{4/(n-2)} \cdot g_M$ and it has \emph{positive mean curvature} $\approx (n-k-1)/\ve$ largely as in the Euclidean case of distance tubes around $\R^k \subset \R^n$. This makes $\p V_{\ve}$ an \emph{outer barrier} for area minimizers in $V_{\ve}$ homologous to $\p V_\ve$. It keeps them inside  $V_\ve$.\\
In turn, for such a (now fixed) $\delta>0$ and $\rho > 0$ small enough we claim that $\p V_\rho$ has \emph{negative mean curvature}. To see this, we
recall that the second fundamental form $A_L(g)$ of a submanifold $L$ with respect to some metric $g$ on the ambient manifold transforms under conformal deformations $u^{4/(n-2)} \cdot g$, for smooth  $u>0$, according to the formula \cite[1.163, p.\,60]{Be}
 \begin{equation}\label{mmm}
 A_L(u^{4/(n-2)} \cdot g)(v,w) = A_L(g)(v,w) - \frac{2}{n-2}\cdot {\cal{N}} (\nabla u / u) \cdot g(v,w),
 \end{equation}
where $ {\cal{N}} (\nabla u / u)$ is the normal component of $\nabla u / u$ with respect to $L$. In our case, we choose $u =\psi$ and locally consider $V_\rho$ as a tube around $\R^k$ in $\R^k \times \R^{n-k}$ where $\R^k$ represents $N$. Thus we have from (\ref{appp}) that $tr\,A_{\p V_\rho}( g ) \approx -(n-k-1)/\rho$, since the $\R^{n-k}$-factor of $V_\rho$ is totally geodesic. The trace of the second summand is  $(n-1) \cdot 2/(n-2)  \cdot (2+k-n)/\rho$ and we get:
\begin{align}
&\qquad\psi^{4/(n-2)} \cdot tr\, A_{\p V_\rho}(\psi^{4/(n-2)} \cdot g)\notag\\
&=  - \left( (n-k-1) +  (n-1) \cdot 2/(n-2)  \cdot (2+k-n)\right) \cdot \rho^{-1}\notag\\
&=(2 - (3 + k) \cdot n + n^2)/(n-2)  \cdot \rho^{-1}\notag\\
&\ge 2/(n-2)  \cdot \rho^{-1}>0, \mm{ since } n-3 \ge k.\label{te}
\end{align}
This makes $\p V_\rho$ an \emph{inner barrier} for area minimizers in $M \setminus V_\rho$ homologous to $\p V_\rho$, keeping them from reaching $N$. Thus for any $\ve >0$, there is some
$\delta_\ve>0$ so that there is a neighborhood $U_{\ve}$ of $N$, with $V_\rho \subset U_{\ve} \subset V_\ve$, such that  $\p U_{\ve}$ is  (locally) area minimizing relative to $g_{\delta_\ve}$.\\

The resulting $scal>0$-manifold $X=X_U$ with its minimal boundary $\p X$ can be used in several ways:\\

\textbf{Doublings}  We  take a second copy $X^*$ of $X$ and glue $X^*$ to $X$ along $\p X$ to get a doubling $X \cup_{\div}  X^*$ with the obvious identification $\div$ along $\p X$. This can be made to get a smooth $scal>0$-metric on $X \cup_{\div}  X^*$ since one can approximate $\p X$, slightly (re)enlarging $X$, by a \textbf{smooth hypersurface $Z$ with positive mean curvature} \cite{G2} and then a standard bending of a collar of $Z$ gives a $scal>0$-metric with totally geodesic boundary $Z$.\\If there is another manifold $Y$ with boundary $\p Y$ being locally isometric to $X^*$ near $\p X^*$ then we can alternatively glue $Y$ to $X$ along $\p X$.
This is just a variant of the classical surgery in \cite{GL} and \cite{SY2}.

\subsubsection{Inductive Removal of Singularities} \label{insw}

As already indicated in the codimension $\ge 3$-surgery reconstructed above, the main result is a partial regularization method where we replace the singularities of $(H^{n},d_{\sima})$ by minimal boundaries:

\begin{theorem}\label{sst}  Let $H^{n}$ be a singular area minimizing hypersurface in some compact $scal \ge 0$-manifold $M^{n+1}$.
Then  there are arbitrarily small neighborhoods $U$ of $\Sigma$, such that $H \setminus U$ is \textbf{conformal} to some $\boldsymbol{scal>0}$\textbf{-manifold} $X_U$  with \textbf{locally area minimizing} boundary $\p X_U$.
\end{theorem}

This is the generic step in an inductive dimensional descent with a built-in partial regularization scheme. Starting from a singular area minimizing hypersurface $H^n$ in some smooth $scal \ge 0$ ambient space  we descent to a possibly again singular area minimizing hypersurface $L^{n-1}$ in the smooth $scal > 0$ ambient space $X_U$. The process shifts the singular issues with  $H^n$ to lower dimensions:\\

\textbf{A.} Let $L^{n-1} \subset \overline{X_U^n}$ be an area minimizer. At this point $\p X_U$ merely constrains $L^{n-1}$ to stay within $\overline{X_U^n}$.
It does not matter whether $\p X_U$ is smooth \cite[Th.1.20, Rm.1.22]{Gi}.

\textbf{B.} The point is that $\p X_U$ is minimal. From this the strict maximum principle  \cite{Si} shows, componentwise,  that either $L^{n-1} \equiv \p X_U^{n-1}$ or $L^{n-1} \cap \p X_U^{n-1} \v$. This renders $L^{n-1}$ as an actually unconstraint minimal hypersurface in an extension $Y_U^n$ of $\overline{X_U^n}$ surrounding $L^{n-1}$ as a smooth $scal>0$-manifold.

\textbf{C.} $L^{n-1}$ may be singular, with a lower dimensional singular set, but we can reapply Theorem \ref{sst} to $L^{n-1}$ in dimension $n-1$. This and further iterations inductively sweep out all singular issues to lower dimensions before they ultimately disappear in dimension $7$.\\

To indicate the proof of this  theorem, we follow the lines of the case without singularities as discussed in the previous section.\\

 \textbf{Step A: Minimal Splitting Factors}\, The counterpart to Step A is the transition from the compact singular area minimizer $H^n$  in some $scal>0$-manifold $M^{n+1}$ to an associated minimal splitting factor  $(H, d_{\sima}(\Phi_H))$, where $\Phi_H>$ is a positive  supersolution of $L_{H,\lambda} \, \phi=0$  for some $0 < \lambda < \lambda^{\bp}_H$. We actually choose a rather small $\lambda$ to get uniform growth estimates for $\Phi$ towards $\Sigma$ which are used to control the secondary deformations in the following Step B. \\

\begin{figure}[htbp]
\centering
\includegraphics[width=1\textwidth]{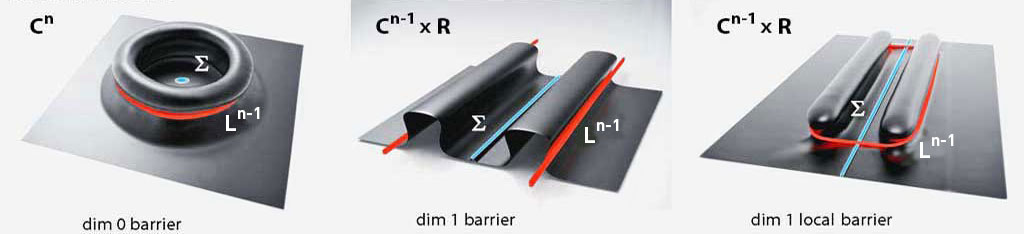}
\caption{The minimal factors are drawn as simple planes. The singular set is visualized as centered spots and lines.  The red rubber bands $L^{n-1}$ illustrate locally area minimizing hypersurfaces kept from contracting to (parts of) $\Sigma$ by these barriers.}
\label{fig:as}
\end{figure}

 \textbf{Step B: Removal of Singularities}\, We deform $(H,d_{\sima})$ in a similar way as in Step B in the smooth sample case above. We find a small neighborhood $U$ of $\Sigma$ and an elementary deformation of $g_{\sima}$ to a $scal >0$-metric on $H \setminus U$ so that $\p U$ becomes \textbf{minimal}.\\
We construct this elementary deformation of $(H,d_{\sima})$ from the top-down analysis of the previous chapter.  We apply tangent cone reductions in the class of minimal splitting factors.  To this end, we cover $\Sigma \subset H$ with  upper bounded intersection numbers by finitely many small balls $B_{r_i}(p_i) \subset H$ which, after scaling to unit size, are well-approximated by the ball $B_1(0) \cap C_i$ in some tangent cone $C_i$ in $p_i$ with its minimal splitting factor geometry.\\
Similarly, we choose a ball cover for $\p B_1(0) \cap \Sigma_{C_i} \subset \p B_1(0) \cap C_i$ and finitely many locally approximating tangent cones. The iteration of this process defines a \textbf{blow-up tree} $\T$ of cones with root $H$. The branching of $\T$ ends, after at most $n-7$ blow-ups, with some product cone $\R^k \times C^{n-k}$ singular only along $\R^k \times \{0\}$.\\
In this case, the minimal splitting factor metric has a simple form. For some function $c:\p B_1(0) \cap C  \ra \R^{>0}$ and a $\gamma \in (-(n-2)/4,0)$, we have \[g_{\sima} = (c(\zeta) \cdot \rho^\gamma)^{4/(n-2)} \cdot g_C \mm{ on } \R^k \times C^{n-k} \setminus \R^k \times \{0\},\] where $\rho(x)=dist(x,\R^k \times \{0\})$, $\zeta(x)=\pi_2(x/|x|)$ and  $\pi_2:\R^k \times C^{n-k}\ra C^{n-k}$ is the projection. For this metric, we easily find $\R^k$-translation invariant conformal deformations, supported away from $\R^k \times \{0\}$, to another $scal>0$-metric which contains an elementary \textbf{inner barrier} for area minimizers illustrated in the left and middle part of \autoref{fig:as}.
An inner barrier is a bumpy deformation of $g_{\sima}$ that locally increases the volume element of  $g_{\sima}$ to keep local area minimizers away from $\Sigma$. These barriers result from a simple cut-off construction of a secondary deformation added to the minimal factor metric.\\
To assemble a \textbf{global barrier} along $\Sigma \subset H$, we localize the barriers in the terminal nodes of $\T$  keeping $scal>0$: we truncate the bump deformation to a deformation supported in a compact subset of
$C \setminus \Sigma_C$, indicated on the right portion of the picture. This particular localization is important for the transfer between different nodes of the blow-up tree. The transfer uses \textbf{smooth} tangent cone approximations, that is, approximations described as sections of the normal bundle of the tangent cone in $C^k$-norm, for some $k \ge 2$, and we use these sections to transfer constructions between different spaces via pull-back or push-forward.  But these smooth approximations exist only in positively upper and lower bounded distance to the singular sets.\\
These transfer maps between nodes in the blow-up  tree allow us to bring the harvest home:  we transport the individual bump deformations backwards to the next node in $\T$
and inductively assemble barriers until we reach $H$. This is a simple process but it needs a considerable amount of bookkeeping using families of well-controlled coverings.
The isoperimetry of $g_{\sima}$ is applied when these deformations have been completed. The  isoperimetric inequality appears in the guise of the volume estimates in Cor.~\ref{horn}.\\

\textbf{Auto-Aligned Coverings:}\, As an alternative to the use of blow-up trees there is a refined
covering argument, in \cite{L5}, that guides the placement of the local barriers.

\begin{itemize}[leftmargin=*]
  \item The  \textbf{local barriers} are placed disjointedly to keep control over their curving effect.  The resulting \textbf{global barrier} along $\Sigma$ does not topologically separate the singular region from the main regular part of $H$. However, using the isoperimetric inequality of the underlying minimal splitting factor, it separates geometrically in the sense that it keeps area minimizers from approaching $\Sigma$. Intuitively, one may compare these topologically permeable barriers with a \textbf{Faraday cage}.
  \item In turn, these volume estimates also ensure that there is such an area minimizer since they show that there is a larger neighborhood of the deformed region that area minimizers do not leave.
\end{itemize}

\textbf{Doublings} \, As in the sample case of codimension $\ge 3$ surgeries, we have arbitrarily small neighborhoods $U$ of $\Sigma$, such that after some smoothing of $\p X_U$ we can take a second copy $X_U^*$ of $X_U$ and glue $X_U^*$ to $X_U$ so that  $X_U \cup_{\div}  X_U^*$ is smooth with $scal>0$.

\subsubsection{From Ordinary to Locally Finite Homology} \label{bom}

Here we will see how to apply the partial regularization method in scalar curvature geometry and general relativity to extend results valid in dimensions $\le 7$ to higher dimensions.
We start with one of the most prominent examples: the Riemannian positive mass theorem. We refer to \cite{L6} for the explicit statement and its physical relevance.
In geometric terms it can be formulated as a \textbf{non-existence of $\boldsymbol{scal > 0}$-islands}:
\begin{theorem}\label{nex} There exists \textbf{no} complete Riemannian manifold $(M^{n+1}, g)$, $n \ge 2$, such that:
\begin{itemize}[leftmargin=*]
   \item $scal(g) >0$ on a non-empty open set $U \subset M^{n+1}$ with compact closure.
   \item $(M^{n+1} \setminus U, g)$ is isometric to $(\R^{n+1} \setminus  B_1(0), g_{Eucl})$.
\end{itemize}
\end{theorem}
We indicate the contradiction argument of \cite{L6}. Let us assume we had such a Riemannian manifold $(M^{n+1}, g)$. Then we can place a large cube around $B_1(0)$ and compactify $(M^{n+1}, g)$
to a closed $scal > 0$-manifold diffeomorphic to $N^{n+1} \cs T^{n+1}$, for some closed manifold $N^{n+1}$. Now we can find a possibly singular area minimizer $H^n \subset N^{n+1} \cs T^{n+1}$ homologous to the standard $T^n \subset T^{n+1}$ and
apply Theorem \ref{sst} to get the $scal>0$-manifold $X_U$ for some small neighborhood $U$ of $\Sigma_H$.\\
The regularity theory also shows that $H^n$ again contains a nearly flat torus summand, even after the conformal deformation of \ref{sst}, i.e., $T^n \setminus B \subset H^n$ for some small ball $B$.  This means that there is a compact area minimizer $L^{n-1}  \subset X_U \subset  H^n$ homologous to $T^{n-1}$ with $\p U \cap L^{n-1} \v$ and, again, $L^{n-1}$ contains a nearly flat torus summand. We iterate this process until we reach a $scal > 0$-surface $F^2 \cs T^2$ that does not exist, due to the Gauss--Bonnet theorem. This proves the non-existence of $scal > 0$-islands and, hence, the positive mass theorem.\\

\textbf{Locally Finite Homology} \, The argument above exploited the existence of a nearly flat torus component of  $N^{n+1} \cs T^{n+1}$ to iteratively bypass the singular sets from a suitable selection of homology classes.

\begin{figure}[htbp]
\centering
\includegraphics[width=1\textwidth]{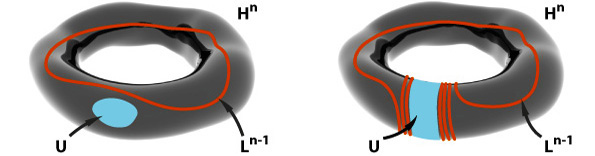}
\caption{The left hand picture represents the case of a compact $L^{n-1} \cap \p U \v$. In the second picture we have a homology class that hits $\overline{U}$. Here we can still find an intrinsically complete local area minimizer $L^{n-1} \cap  \p U \v$ that asymptotically approaches $\p U$.}
\label{fig:to}
\end{figure}

Theorem \ref{sst} can also be used in a different way to include more general homology classes. We start with any non-trivial family of classes $\alpha[1],\dots,\alpha[k] \in H^1(M^{n+1},\Z)$. Geometric measure theory shows that  there is an area minimizer $H^n \subset M^{n+1}$ that represents $\alpha[1]  \cap [M]$ in the homology of integral currents. Since $M$ is smooth, the integral current homology groups of $M$ are isomorphic to ordinary homology groups \cite{D}, but, in general, a singular $H^n$ does \emph{not} represent the associated class in ordinary homology.
\begin{itemize}[leftmargin=*]
  \item $H^n$ might not be a singular $n$-cycle or have a non-finitely generated homology and $H^n \setminus \Sigma$ does not represent a class in ordinary homology.
\end{itemize}
This advocates the interpretation of the $\alpha[1] \cap \cdots \cap \alpha[m] \cap [M]$, $m \le k$, as classes in integral current homology. However, Theorem \ref{sst} and the regularity theory for minimal hypersurfaces, proving the low dimensionality of the singular set, suggest still another interpretation. We think of the $\alpha[1] \cap \cdots \cap \alpha[m] \cap [M]$ as classes in a \emph{locally finite} homology \cite{HR}. \\
In more detail, for small $U$, the isoperimetric inequality and the fact that $\p U$ is locally area minimizing can be used to show that there is an intrinsically complete local area minimizer  $L^{n-1} \subset Y_U \subset H^n$ so that, componentwise, either $L^{n-1}$ is compact and $L^{n-1} \cap \p X_U \v$ (or $L^{n-1} =\p X_U$), this is the case we have already discussed above, or
\begin{itemize}[leftmargin=*]
\item  $L^{n-1}$ is  \textbf{non-compact with end}. It asymptotically approaches $\p X_U$, like a geodesic ray on a surface that approaches a closed geodesic in infinitely many loops. We apply Theorem \ref{sst}  to $\p X_U$ and eventually transfer its metric to  $L^{n-1}$ to get a periodic $scal>0$-end structure.
\end{itemize}
For the present we assume that $L^{n-1}$ is smooth. Then $L^{n-1}$ represents a class in the Borel--Moore homology $H_{n-1}^{BM}(X_U)$ of $X_U$. $L^{n-1}$ can be thought to represent $\alpha[1]  \cap \alpha[2] \cap [M]$. To make this intuition more precise, the local finiteness of chains in this homology becomes  important. $L^{n-1}$ spins around $\p X_U$, but each of the infinitely many loops can be annihilated by one ordinary homology operation. That is, $L^{n-1}$  is Borel--Moore homologous to the restriction $L_*^{n-1} \cap X_U$  of a hypersurface $L_*^{n-1} \subset H$ representing  $\alpha[1]  \cap \alpha[2] \cap [M]$ in integral current homology.\\
  In our context, this suffices to view $L^{n-1}$ as a representing cycle of $\alpha[1]  \cap \alpha[2] \cap [M]$  since area minimizers reaching into the periodic end of $L^{n-1}$ can be made to stay in a compact subset of $L^{n-1}$. A sample of such a truncation style process is explained in the proof that topologically large manifolds cannot admit $scal>0$-metrics \cite{L7}. \\
  \emph{It would be nice to formalize such ad-hoc truncation arguments in a suitable locally finite or coarse homology theory where $L^{n-1}$ properly represents $\alpha[1]  \cap \alpha[2] \cap [M]$, cf.\ \cite{HR}, \cite{R} for some background.}\\
To finish the argument, we inductively apply the same strategy to get potentially singular and non-compact area minimizers with ends that, in the sense above, represent $\alpha[1]  \cap\cdots\cap \alpha[m] \cap [M]$ in smooth $scal>0$-manifolds of dimension $n-m+1$. Then we regularize them, using Theorem  \ref{sst}, before we start the next loop. This way, we inductively sweep out all singular issues to lower dimensions before they ultimately disappear in dimension $7$.

\end{document}